\documentclass[12pt]{article}
\usepackage{amsxtra,amssymb,amsthm,amsmath,latexsym}

\theoremstyle{plain}
\newtheorem{theorem}{Theorem}

\newtheorem{lemma}{Lemma}%[section]
\newcommand{\refT}[1]{Theorem~\ref{T:#1}}
\newcommand{\refS}[1]{Section~\ref{S:#1}}

\newcommand{\refL}[1]{Lemma~\ref{L:#1}}

\def\dotu{\dot{u}}

\def\oH{{\overset{\circ}{H}}}
\def\oH1{{\overset{\circ}{H}\kern-.02in{}^1}}

\def\ve{{\varepsilon}}
\def\bee{\begin{equation*}}
\def\eee{\end{equation*}}
\def\be{\begin{equation}}
\def\ee{\end{equation}}

\begin{document}
%\begin{titlepage}
\title{                  
Iterative solution of linear equations with unbounded operators}

\author{
A.G. Ramm\\
 Mathematics Department, Kansas State University, \\
 Manhattan, KS 66506-2602, USA\\
ramm@math.ksu.edu\\
http://www.math.ksu.edu/\,$\widetilde{\ }$\,ramm}

\date{}

\maketitle\thispagestyle{empty}

\begin{abstract} 
\footnote{Math subject classification: 44A05 \quad 47A50 \quad 65J10}
\footnote{key words:  iterative methods, unbounded operators, linear
equations.}

A convergent iterative process is constructed for solving any solvable
linear equation in a Hilbert space.
\end{abstract}

%\end{titlepage}

\section{Introduction}\label{S:1}

A basic general result about solvable linear equations
\be\label{e1} Au=f, \ee
where $A$ is a linear bounded operator in a Hilbert space, is the following
theorem.

\setcounter{theorem}{-1}
\begin{theorem}\label{T:0}
Any solvable equation \eqref{e1} with a bounded linear operator can be 
solved by a convergent iterative process.
\end{theorem}

A proof of \refT{0} can be found, e.g., in \cite{2}.
One of the steps in this proof is the following simple Lemma
(see e.g. \cite{2}).

\begin{lemma}\label{L:1}
If equation \eqref{e1} is solvable and $A$ is a bounded linear operator, 
then equation \eqref{e1} is equivalent to
\be\label{e2} A^\ast Au=A^\ast f. \ee
\end{lemma}

The equivalence in \refL{1} means that every solution to
\eqref{e1} solves \eqref{e2} and vice versa.

The aim  of this paper is to study equation \eqref{e1} with a linear,
closed, densely defined, unbounded, and not necessarily boundedly
invertible operator. In other words, $A$ may be not injective, i.e.,
null space $N:=N(A)$ may be non-trivial, and its range $R(A)$ may be
not closed. Although there are many papers and books on iterative methods,
iterative methods for equations \eqref{e1} with unbounded operators were not
studied in such generality.

Our second aim is to study a variational regularization method for the
solutions to equation \eqref{e1}. By $y$ we denote throughout the unique
solution to \eqref{e1} of minimal norm, i.e., the solution $y\perp N$.
This solution will be of main interest to us. If $A$ is bounded, but
not boundedly invertible, so that \eqref{e1} is an ill-posed problem
(see e.g. \cite{2}), then a variational regularization method for obtaining a
stable approximation of the solution $y$ given noisy data $f_\delta$,
$\|f_\delta-f\|\leq\delta$, consists of

a) minimizing the functional
\be\label{e3} F(u)=\|Au-f_\delta\|^2+a\|u\|^2, \qquad a>0, \ee
where $a$ is a constant called a regularization parameter, proving that
\eqref{e3} has a unique global minimizer
$u_{a,\delta}=(A^\ast A+aI)^{-1}A^\ast f_\delta$, 

and

b) proving that one can choose $a=a(\delta)$, so that
$\lim_{\delta\to 0}a(\delta)=0$ and
\be\label{e4}
  \lim_{\delta\to 0}\|u_\delta-y\|=0,
  \qquad u_\delta:=u_{a(\delta),\delta}.\ee

Formula \eqref{e4} shows that $u_\delta$ is a stable approximation of $y$.
The rate of convergence of $u_\delta$ to $y$ is not possible to specify
without imposing additional assumptions on $f$.

If $A$ is unbounded, then it was not proved that functional \eqref{e3} has
a unique global minimizer. Formula
$u_{a,\delta}=(A^\ast A+aI)^{-1} A^\ast f_\delta$ 
is not well defined
because $f_\delta$ may not belong to $D(A^\ast)$.

Throughout the paper $T=A^\ast A$ is a selfadjoint nonnegative
operator (generated by the closed nonnegative quadratic form
$(Au,Au)$, $D(T)\subset D(A)$, $T_a:= T+aI$, $I$ is the identity operator,
$Q=AA^\ast\geq 0$ is a selfadjoint operator, $D(Q)\subset D(A^\ast)$.
Recall that $A^\ast$ is well defined if $A$ is densely defined,
and $A^\ast$ is densely defined if $A$ is closed. (See \cite{1}).
By $S_0$ we denote the operator $T^{-1}_{a}A^\ast$ with domain 
$D(A^\ast)$, and by $S$ we denote its closure.

Our results can be described as follows: the operator $S_0$ is closable,
its closure is defined on all of $H$ and is a bounded operator,
$\|S\|\leq\frac{1}{2\sqrt{a}}$.
A similar result holds for $S_1=T^{-1}_{ia} A^\ast$ with domain 
$D(A^\ast)$.
Our result shows that the element $T^{-1}_a A^\ast f_\delta$ is well
defined for any $f_\delta\in H$, and not only for $f_\delta\in D(A^\ast)$.

Consider the iterative process
\be\label{e5}
  u_{n+1}=Bu_n+T^{-1}_a A^\ast f, \quad u_1=u_1, \quad a>0, \ee
where the initial approximation $u_1\perp N$ and otherwise arbitrary, and 
$B:=a T^{-1}_a$.

\begin{theorem}\label{T:1}
If $A$ is a linear, closed, densely defined operator in $H$, $a>0$, 
$B=aT^{-1}_a$, and $u_1\perp N$, then
\be\label{e6} \lim_{n\to\infty}\|u_n-y\|=0,\ee
where $u_n$ is defined by (5).
\end{theorem}

\begin{theorem}\label{T:2}
If $A$ is a linear, closed, densely defined operator in $H$, and $a>0$ is 
a constant, then the operator $S_0=T^{-1}_{a} A^\ast$ with domain 
$D(A^\ast)$ is closable and its closure $S$ is a bounded operator defined 
on all of $H$, $\|T^{-1}_a A^\ast\|\leq\frac{1}{2\sqrt{a}}$. Similar 
results hold for the operator $T^{-1}_{ia} A^\ast$.
\end{theorem}

\begin{theorem}\label{T:3}
For any $f_\delta\in H$ functional \eqref{e3} has a unique global 
minimizer
$u_{a,\delta}=T^{-1}_a A^\ast f_\delta=A^\ast Q^{-1}_a f_\delta$.
\end{theorem}

In Section 2 proofs are given. In Section 3 we construct a
stable approximation to $y$ given noisy data $f_\delta$ and
using an iterative process similar to \eqref{e5}.
In Section 4 the case of selfadjoint, unbounded and possibly 
not boundedly invertible operator is briefly considered. In 
Section 5 the dynamical systems method (DSM) (developed in 
[2] pp.41-70) is justified for equation (1) with unbounded, 
linear, densely defined operator in a Hilbert space.
The basic results of this paper are stated in Theorems 1 
through 5.

\section{Proofs} \label{S:2}

\begin{proof}[Proof of \refT{2}]
To prove the closability of $S_0$, assume that $u_n\in D(A^\ast)=D(S_0)$,
$u_n\to 0$, $S_0 u_n\to f$, and prove $f=0$. We have
\be\label{e7}
  (f,h)=\lim_{n\to\infty}(S_0 u_n,h)
   =\lim_{n\to\infty} (A^\ast u_n, T^{-1}_a h)
   =\lim_{n\to\infty}(u_n,AT^{-1}_ah)=0,
   \qquad \forall h\in H. \ee
Here we have used the inclusion $R(T^{-1}_a)\subset D(A)$.
This inclusion can be verified: if $g=T^{-1}_ah$, then
$Tg+ag=h$, so $g\in D(T)\subset D(A)$, as claimed.
From \eqref{e7} it follows that $f=0$ because $h\in H$ is arbitrary.
Relation \eqref{e7} shows that $D(S^\ast_0)=H$ and $S^\ast_0=AT^{-1}_a$.
This operator is closed and densely defined. Indeed,
by the polar decomposition, $A=UT^{\frac{1}{2}}$, where $U$ is a partial
isometry, $\|U\|\leq 1$. The operator $T^{\frac{1}{2}} T^{-1}_{a}$
is densely defined, it is a function of the selfadjoint operator $T$.
We have
\be\label{e8}
  \|S^\ast_0\| \leq \|T^{\frac{1}{2}} T^{-1}_a\|
  = \| \int^\infty_0 \frac{s^{\frac{1}{2}}}{s+a} dE_s \|
  =\sup_{s\geq 0} \frac{s^{\frac{1}{2}}}{s+a}
  =\frac{1}{2\sqrt{a}}, \ee
where we have used the spectral theorem and $E_s$ is the resolution
of the identity of the selfadjoint operator $T$.
Since $\|S\|=\|S^{\ast\ast}_0\| = \|S^\ast_0\| \leq \frac{1}{2\sqrt{a}}$,
\refT{2} is proved except for the claim concerning the operator
$T^{-1}_a A^\ast$. The proof of this claim is essentially the same as
the above proof,
the only (not important) difference is the replacement of
the formula $(T^{-1}_a)^\ast=T_a^{-1}$ by 
$(T^{-1}_{+ia})^\ast=T^{-1}_{-ia}$.

\refT{2} is proved.
\end{proof}

\begin{proof}[Proof of \refT{1}]
If equation \eqref{e1} is solvable then $f=Ay$, the operator $B=a 
T_{a}^{-1}$
is bounded, defined on all of $H$, and $\|B\|\leq 1$.

Consider the equation
\be\label{e9}  u=Bu+T^{-1}_a A^\ast f.\ee
This equation makes sense for any $f\in H$ by \refT{2}.
The minimal-norm solution $y$ to equation \eqref{e1} solves \eqref{e9} in
the following sense.

If $f\in D(A^\ast)$, then $Ty=A^\ast f$, $ay+Ty=ay+A^\ast f$,
$y=By+T^{-1}_a A^\ast f$, so $y$ solves \eqref{e9} if $f\in D(A^\ast)$.
Since the set $f\in D(A^\ast) \bigcap R(A)=D(T)$ is dense in $H$ and,
consequently, in $D(A^\ast)$, and since, by \refT{2}, the operator
$S=T^{-1}_a A^\ast$ is uniquely extendable to all of $H$
(from a dense subset $D(A^\ast)$) by continuity, it follows that if
$y$ solves equation \eqref{e9} for every $f\in D(A^\ast)\bigcap R(A)$,
then this equation is solvable for any $f\in R(A)$.
Indeed, suppose $u_n= Bu_n+Sf_n$, $f_n\in D(A)\bigcap R(A)$, $u_n\perp N$,
and $\lim_{n\to\infty} f_n=f\in R(A)$. The subspace
$N^\perp:=N(T)^\perp =\overline{R(T)}$ is invariant with respect to $B$.
If $f\in R(A)$ and $f=Ay$, then equation \eqref{e9} has a unique solution
$u\in N^\perp$, and this solution is $u=(I-B)^{-1} T^{-1}_a Ty=y$.

Indeed,
\bee
  (I-B)^{-1} T^{-1}_{a} Ty=(I-aT^{-1}_a)^{-1} T^{-1}_a Ty
  =\int^\infty_0 \frac{s dE_sy}{(1-\frac{a}{a+s})(s+a)}=y \eee
Denote $w_n:=u_n-y$. Then \eqref{e5} and equation \eqref{e9} for $y$
imply $w_{n+1}=Bw_n$, so $w_{n+1}=B^nw$, $w:=u_1-y$, $w\perp N$.

Let us prove $\lim_{n\to\infty}\|B^n w\|=0$.
If this is proved, then \eqref{e6} follows,
and \refT{1} is proved. We have
\be\label{e10}
 I:=\|B^nw\|^2 = \int^\infty_{0} \frac{a^{2n}}{(a+s)^{2n}} d(E_sw,w)
 =\int_{s\geq b} + \int_{0<s\leq b}:=I_1+I_2 \ee
where $b>0$ is a number and $E_s$ is the resolution of the identity
corresponding to the selfadjoint operator $T\geq 0$.

In the region $|s|\geq b$ one has
\be\label{e11}
  I_1\leq q^n(b) \|w\|^2, \qquad 0<q<1,
  \qquad q=q(b)=\max_{|s|\geq b} \frac{a^2}{(a+s)^2}.\ee
We estimate $I_2$ as follows:
\be\label{e12}
  I_2\leq \int^b_{-b} d(E_sw,w) \ee
uniformly with respect to $n$.

Since $w\perp N$, we have
\be\label{e13}
  \lim_{b\to 0} \int^b_0 d(E_sw,w)=\|P_Nw\|^2=0, \ee
where $P_N=E_{+0}-E_0$ is the orthoprojector onto $N$, and
\be\label{e14}
  \lim_{b\to 0} \int^0_{-b} d(E_sw,w)=0 \ee
because $E_{s-0}=E_s$. Therefore, for an arbitrary small $\ve>0$,
we choose $b>0$ so small that $I_1\leq\frac{\ve}{2}$, and for fixed $b$
we choose $n$ so large that $q^n(b)\leq\frac{\ve}{2}$.
Then $I\leq\ve$.

\refT{1} is proved.
\end{proof}

We could replace $a$ by $ia$ in the above arguments.

\begin{proof}[Proof of \refT{3}]

Denote $f_\delta:=g$ and $u_{a,\delta}:=z:=A^\ast Q^{-1}_a g$.
The operator $A^\ast Q^{-1}_a=V Q^{\frac{1}{2}} Q^{-1}_a$, where $V$
is a partial isometry, so
$\|A^\ast Q^{-1}_a\|\leq \|Q^{\frac{1}{2}} Q^{-1}_a\|
\leq \frac{1}{2\sqrt{a}}$.
Thus, $z$ is defined for any $g\in H$. We have
\be\label{e15}
  F(z+h)=F(z)+ \|Ah\|^2 + a\|h\|^2 + 2Re[(Az-g,Ah)+a(z,h)]
  \qquad \forall h\in D(A).\ee
If $z=A^\ast Q^{-1}_a g$, then
\be\label{e16}
  (Az-g,Ah) + a(z,h)=(QQ^{-1}_a g-g,Ah)+a(Q^{-1}_a g,Ah)=0.\ee
From \eqref{e15} and \eqref{e16} we obtain
\be\label{e17}
  F(z+h)\geq F(z) \qquad \forall h\in  D(A) \ee
and $F(z+h)=F(z)$ implies $h=0$.
Thus $z$ is the unique global minimizer of $F$.
Let us prove $A^\ast Q^{-1}_a =T^{-1}_a A^\ast$. Since both operators in
this identity are bounded, it is sufficient to check that
\be\label{e18}
  A^\ast Q^{-1}_a \psi=T^{-1}_a A^\ast \psi \ee
for all $\psi$ in a dense subset of $H$. As such dense subset let us take
$D(A^\ast)$. Denote $Q^{-1}_a \psi:=g$. Then $\psi=Q_a g$.
Equation \eqref{e18} is equivalent to
\bee
  T_aA^\ast g=A^\ast Q_a g, \hbox{\ or\ }
  A^\ast A A^\ast g+aA^\ast g=A^\ast AA^\ast g+aA^\ast g,\eee
which is an obvious identity. Reversing the steps, we obtain \eqref{e18}
for every $\psi\in D(A^\ast)$. Note that $\psi\in D(A^\ast)$ is equivalent
to $g\in D(A^\ast AA^\ast)$, so that the above calculations are justified.

\refT{3} is proved.
\end{proof}

Equation \eqref{e18} allows one to replace the term
$T^{-1}_a A^\ast f$ in \eqref{e5} by the term $A^\ast Q^{-1}_a f$
which is originally well defined for any $f\in H$.

\section{Stable solution of (1).} \label{S:3}
Suppose that noisy data $f_\delta$, $\|f_\delta-f\|\leq\delta$, are given.
We want to construct a stable approximation $u_\delta$ of $y$ in the sense
\eqref{e4}.

One way to do this is to use iterative process \eqref{e5}, with $f_\delta$
in place of $f$, and stop the iterations at the step $n=n(\delta)$,
where $n(\delta)$ is properly chosen.
Indeed, if we use the argument from the proof of \refT{3}, then we get
$w_{n+1}=Bw_n+S(f_\delta-f)$.
Thus $w_{n+1}=\sum^n_{j=0} B^j S(f_\delta-f) + B^nw$, where
$w=u_1-y$, $u_1\perp N$, and
$$\|\sum^n_{j=0} B^j S(f_\delta-f)\|
 \leq \frac{\delta(n+1)}{2\sqrt{a}} + \ve(n):=\nu(\delta,n)$$
where $\lim_{n\to\infty}\ve(n)=0$, as we have demonstrated in
the proof of \refT{1}.

It is clear that if $n=n(\delta)$ is chosen so that
$\lim_{\delta\to 0}n(\delta)=\infty$ and
$\lim_{\delta\to 0}\delta n(\delta)=0$,
then $\lim_{\delta\to 0}\nu(\delta,n(\delta))=0$.
Therefore $u_\delta=u_{n(\delta)}$ satisfies \eqref{e4}.

\section{Equation (1) with selfadjoint operator.}\label{S:4}
Assume that $A=A^\ast$, $A$ is unbounded and $A$ does not have a
bounded inverse. Then one can use an analog of \refT{1} in the form
\be\label{e19}
  u_{n+1}=L u_n+g, \qquad L:=ia(A+ia)^{-1}, \qquad g=(A+ia)^{-1} f,\ee
where $a=const>0$ and $u_1\perp N:=N(A)$ is arbitrary. Note that any 
element
$Ah$, $\forall h\in D(A)$, is orthogonal to $N$ since $A=A^\ast$.
The minimal-norm solution $y$ to \eqref{e1} solves the equation
$y=Ly+g$, so that the proof of \refT{1} remains almost the same.
So we get

\begin{theorem}\label{T:4}
If $A=A^\ast$ is unbounded, $a=const>0$, and equation \eqref{e1} is solvable,
then \eqref{e6} holds for the iterative process \eqref{e19}.
\end{theorem}

Also an analog of the result of \refS{3} holds.

\section{DSM}\label{S:5}
In this Section we justify the dynamical systems method (DSM)
for solving equation \eqref{e1}. The DSM theory is developed in [2], 
pp.41-70.

\begin{theorem}\label{T:5}
Assume that $f=Ay$, $y\perp N$, $N:=N(A)$, $A$ is a linear operator, 
closed and densely
defined in $H$. Consider the problem
\be\label{e20}
  \dotu= -u+T^{-1}_{\ve(t)} A^\ast f, \qquad u(0)=u_0;\qquad 
\dotu:=\frac {du}{dt},\ee
where $u_0\in H$ is arbitrary, $T_\ve=A^\ast A+\ve I$,
$\ve=\ve(t)>0$ is a continuous function monotonically decaying to zero
as $t\to\infty$, and $\int^\infty_0\ve(s)ds=\infty$.
Then problem \eqref{e20} has a unique solution $u(t)$
defined on $[0,\infty)$, there exists
\be\label{e21}
  \lim_{t\to\infty}u(t):=u(\infty), \quad \hbox { and }\quad u(\infty)=y. 
\ee
\end{theorem}

\begin{proof}
One has  $T^{-1}_{\ve(s)}A^*f=T^{-1}_{\ve(s)}Ty$. Therefore
\be\label{e22}
  u(t)=u_0 e^{-t}+\int^t_0 e^{-(t-s)} T^{-1}_{\ve(s)} Ty\,ds.\ee

The conclusion of \refT{5} follows immediately from two lemmas:
\setcounter{lemma}{1}
\begin{lemma}\label{L:2}
If there exists $h(\infty)=\lim_{t\to\infty} h(t)$, then
\be\label{e23}
  \lim_{t\to\infty} \int^t_0 e^{-(t-s)} h(s)\,ds=h(\infty).\ee
\end{lemma}

\begin{lemma}\label{L:3}
If $y\perp N:=N(A)$, then
\be\label{e24}
  \lim_{\ve\to 0} T^{-1}_\ve Ty=y.\ee
\end{lemma}

The proof of \refL{2} is simple and  is left to the reader.

The proof of \refL{3} is briefly sketched below::
\bee
  T^{-1}_\ve Ty-y= \int^\infty_0\left(\frac{s}{s+\ve}-1\right)
  d E_s y= -\int^\infty_0 \frac{\ve}{s+\ve} d E_s y.\eee
Thus,
\be\label{e25}
  \lim_{\ve\to 0} \|T^{-1}_\ve Ty-y\|^2
  = \lim_{\ve\to 0} \int^\infty_0 \frac{\ve^2d(E_sy,y)}{(s+\ve)^2}
  =\|P_Ny\|^2=0, \ee
because $y\perp N$.
The projector $P_N$ is the orthogonal projector onto $N$.

\refT{5} is proved.
\end{proof}

One can use \refT{5}, exactly as it is done in \cite{2}, for stable
solution of equation \eqref{e1} with noisy data: if $f_\delta$ is
given in place of the exact data $f$, $\|f_\delta-f\|\leq\delta$,
then one solves problem \eqref{e20} with $f_\delta$ in place of $f$,
calculates its solution $u_\delta(t)$ at $t=t_\delta$, and proves that
\be\label{e26}
  \lim_{\delta\to 0} \|u_\delta(t_\delta)-y\|=0 \ee
if $t_\delta$ is suitably chosen. The stopping time $t_\delta$ can be
uniquely determined, for example, by a discrepancy principle
as shown in \cite{2} for bounded operators $A$.
The argument in \cite{2} remains valid in the case of unbounded $A$
without any changes.


\begin{thebibliography}{1000} %number of characters of longest bibitem label

%\bibitem{R1} autonumbers \cite{R1}
%\bibitem[R1]{R1} prints [R1] for \cite{R1}

\bibitem{1}
Kato, T.,
\textbf{Perturbation Theory for Linear Operators},
Springer Verlag, New York, 1984.

\bibitem{2}
Ramm, A.~G.,
\textbf{Inverse Problems}, Springer Verlag, New York, 2005.

\end{thebibliography}
\end{document}